\documentclass[a4paper,12pt]{amsart}
\usepackage{amsmath,amsthm,amssymb,hyperref,graphicx,mathrsfs,mathtools}  
\theoremstyle{plain}
\newtheorem{theorem}{Theorem}

\allowdisplaybreaks
\numberwithin{equation}{section}
\title[Sendov's Conjecture]{A note on a recent attempt to prove Sendov's conjecture}
\author{N.A. Rather \& Suhail Gulzar}
\address{Department of Mathematics, University of Kashmir, Srinagar-190006, India}
\email{dr.narather@gmail.com, sgmattoo@gmail.com}
\begin{document}
\maketitle
\begin{abstract}
Recently G.M. Sofi \& S.A. Ahangar \cite{ss} made an attempt to prove the Sendov's conjecture. But unfortunately the proof is not correct. In this note, we discuss the fallacy in the proof.
 \end{abstract}
Sendov's conjecture says that if all roots of a polynomial $p(z)$ lie within the unit disk, then for every root, there exists a critical point at a distance at most one from the root. Since the Gauss-Lucas theorem implies that the critical points of $p(z)$ must themselves lie in the unit disk, it seems completely implausible that the conjecture could be false. Yet, at present, it has not been proven for polynomials with real coefficients or for any polynomial whose degree exceeds $8.$

Recently G.M. Sofi \& S.A. Shabir \cite{ss} claimed to have proved Sendov's conjecture. But unfortunately the proof is not correct. The proof is divided into two cases. In the first case they consider the class of monic polynomial $p(z)$ such that $\max_{|z|=1}|P(z)|\leq 1.$ The main flaw in the proof lies in this hypothesis. It is easily to observe that if $p(z)$ is a monic polynomial then ${\max}_{|z|=1}|P(z)|\geq 1$ and that equality holds if and only if $p(z)=z^n.$ To see this:

Consider $p(z)=z^n+a_{n-1}z^{n-1}+\cdots+a_0$ and $q(z)=z^np(1/z)$ then
\begin{align*}
\max\limits_{|z|=1}|p(z)|=\max\limits_{|z|=1}|q(z)|\geq |q(0)|=1.
\end{align*}
Now, if $p(z)=z^n$ then clearly $\max_{|z|=1}|p(z)|=1$. Suppose that $\max_{|z|=1}|p(z)|=1$, then also $\max_{|z|=1}|q(z)|=1$. By the maximum principle, we have that $q(z)\equiv 1 $ and this implies $p(z)=z^n$. 

This fact can also be observed by the following \cite{sg} generalization of Visser's inequality \cite{v}.
\begin{theorem}
If $p(z)=\sum_{j=0}^{n}a_jz^j$ be a polynomial of degree $n,$ then 
\begin{align*}
\underset{|z|=1}{\max}|p(z)|\geq |a_n|+\frac{|a_k|}{\binom{n}{k}},\quad k=0,1,\ldots, n-1.
\end{align*}
\end{theorem}
This clearly shows that if $p(z)$ is a monic polynomial such that $p(z)\neq z^n$ then ${\max}_{|z|=1}|p(z)|>1.$

In the second case of the proof they consider the case $\max_{|z|=1}|p(z)|>1$ and apply Case-I to $f(z)=\frac{p(z)}{\max_{|z|=1}|p(z)|}$.  Since $f(z)$ is no longer a monic polynomial then the application of Case-I to $f(z)$ is not valid. Thus this case is also incorrect.
 
.

\end{document}